\documentclass{alf-j-l}
\usepackage{verbatim} 
\usepackage{url} 
   \def\K{{\mathbb K}}
\def\Z{{\mathbb Z}} \def\Q{{\mathbb Q}} \def\A{\overline{\mathbb Q}}

\def\F{{\mathbb F}} 

\def\({\bigl(} \def\){\bigr)} 
\def\\{\cr} 
\def\lt{<}

\def\ov{\overline}

\def\CF#1{{\def\\{\mathrel{,}}\def\;{\mathrel{;}}
\def\dots{\ldots}\def\dotss{\ldots\ldots}[\,#1\,]}}
\def\lCF#1{{\def\\{\mathrel{,}}\def\;{\mathrel{;}}
\def\dots{\ldots}\def\dotss{\ldots\ldots}[\,#1\,,}}

\def\divides{{\mathchoice{\mathrel{\bigm|}}{\mathrel{\bigm|}}{\mathrel{|}}%
{\mathrel{|}}}}
\def\Div{\divides}
\def\notdivides{\mathrel{\kern-3pt\not\!\kern3.5pt\bigm|}}

\def\poly {polynomial}

\def\cf {continued fraction}
\def\pq {partial quotient}

\def\ex{expansion}

\def\cfe{continued fraction expansion}

\DeclareMathOperator{\Jac}{Jac}

\theoremstyle{plain}  
\newtheorem{theorem}{Theorem}
\newtheorem{lemma}[theorem]{Lemma}

\newtheorem*{reduction}{Reduction Theorem}
\newtheorem{corollary}[theorem]{Corollary}

\theoremstyle{definition}  
\newtheorem{definition}[theorem]{Definition}

\newtheorem{example*}{Example}

\theoremstyle{remark}

\newtheorem*{remark*}{Remark}

\copyrightinfo{\number\year}{Alfred J van der Poorten}

\begin{document}

\title {Specialisation and
Reduction\\ of 
\\Continued
Fractions of Formal Power Series \\\phantom{nothing}} 

\author{Alfred J. van der Poorten}
\address{\hskip-1.4\parindent ceNTRe for Number Theory Research,
1 Bimbil Place, Killara, Sydney 2071, 
\penalty10000 Australia}
\email{alf@math.mq.edu.au (Alf van der Poorten AM)}

\thanks{The author was supported in part by a grant from the Australian
Research Council.}


\date{\today.}

\dedicatory{To Jean--Louis Nicolas in celebration of his sixtieth birthday}


\begin{abstract} We discuss and illustrate the behaviour of the
\cfe\ of a formal power series under specialisation of parameters or their
reduction modulo~$p$ and sketch some applications of the reduction
theorem here proved.
\end{abstract}

\maketitle
\pagestyle{myheadings}\markboth{{\headlinefont Alf van der
Poorten}}{{\headlinefont Specialisation of Continued Fraction
Expansions}}

\section{Introduction}\noindent
Given a \cfe\ over a function field $\K(X)$
it may happen that a specialisation of incidental parameters leads
some \pq\ to acquire infinite coefficients. I point out that in such
a case the \pq\ will `collapse' to a \pq\ of higher degree. Indeed,
I illustrate several techniques for manipulating \cfe s to display
such a `collapse' explicitly.

To make sense of the notion ``specialisation'' the base field
$\K=K(t_1,t_2,\ldots)$ should be a transcendental extension of some yet
more base field $K$ by one or several algebraically independent parameters
$t_1$, $t_2$, $\ldots\,$.  In the
sequel $K$ will be supposed to be $\Q$ or some finite field $\F_p$, but our
remarks are often more general.

For example, if $Y^2=(X^2+u)^2+4v(X+w)$
viewed as defined over the field $K$, then studying the continued fraction
of $Y$ will have us working over $\K(X)$, where $\K=K(u,v,w)$.

\begin{definition}  A \emph{specialisation} is a restriction
on the generality of the parameters defining $\K$ over
$K$.\end{definition}
\noindent For example, taking $u+w^2=v$ above is a specialisation; so is
a parametrisation $u=1-2t-t^2$, $v=1-2t$, $w=t$. \emph{A fortiori}, a
numerical example is a specialisation.

Our $Y$ above is a formal Laurent series, an element
of $\K((X^{-1}))$. Its `Taylor' coefficients are elements of $\K=K(u,v,w)$. 
For simplicity, let me suppose $K=\Q$. I say that replacing the base field
$\Q$ by some finite field, say $\F_p$ is  a reduction of $Y$. Of course such
a reduction may not make sense, because $p$ occurs in the denominator of
some Taylor coefficient. If so, I say variously that $Y$ does not have
reduction modulo~$p$ or has bad reduction at~$p$. 

Furthermore, one sees readily that for the
example $Y^2=(X^2+u)^2+4v(X+w)$ reduction modulo~$2$ either yields $Y=X^2+u$,
a \poly\ or, for specialisations of $u$, $v$, $w$, may make no sense. In
either case I may call the reduction bad, in the former case because under
reduction the curve decreases in genus.

\begin{definition} A reduction of an element  $Y\in\K((X^{-1}))$ modulo
$\mathfrak p$ is the replacement, if that makes sense, of $Y$ by the
corresponding series defined over $\K_{\mathfrak p}((X^{-1}))$. Here
$\mathfrak p$ is a maximal ideal of the ring of integers $\mathcal O_{\K}$ of
$\K$  and
$\K_{\mathfrak p}=\mathcal O_{\K}/{\mathfrak p}\mathcal O_{\K}$.
\end{definition}

\noindent 

These remarks serve as introduction to the principal result of this
note.

\begin{reduction}\label{eq:mainthm}
Suppose $F=\sum f_iX^{-i}$ is a formal series in $\K((X^{-1}))$, say
with
\cfe\
$\CF{a_0(X)\\ a_1(X)\\\ldots\,}$ and thus with
convergents given for $h=0$, $1$, $\ldots\,$ by the rational
functions
$$r_h(X)=x_h(X)/y_h(X)=\CF{a_0(X)\\ a_1(X)\\\ldots\\a_h(X)}.
$$
Now denote by $\ov F=\sum \ov f_iX^{-i}$ some
specialisation of $F$ effected by
specialising parameters of the field over which $F$ is
defined, or a reduction of $F$. Then the sequence of convergents to $\ov F$
is precisely the sequence $\(\ov r_h(X)\)$, however listed \emph{without}
repetition. \end{reduction}
Note here that if some \pq\ $a_h$ has no specialisation (in more blunt
words, blows up) then neither do the
\emph{continuants} $x_h$ and $y_h$. Their quotients
$r_h$ will nevertheless have a defined specialisation provided only that the
specialisation, or reduction, $\ov F$ makes sense. In practice, several
distinct convergents
$r_h$, $r_{h+1}$, $\ldots\,$, may collapse to the same convergent
$\ov r_h$ of
$\ov F$.

A special case of this phenomenon is noted in \cite{Ca}. Here there is
no collapse on reduction modulo~$3$. Cantor remarks that because all
the partial quotients of 
$$G_3(X)=\prod_{h=0}^\infty\(1+X^{-3^h}\)
$$
are linear if $G_3$ is viewed as defined over $\F_3$ they must also
be linear for $G_3$ defined over $\Q$. Specifically, in that example
the \pq s of $G_3$ over $\Q$ all have good reduction modulo~$3$ (as
was  shown with inappropriate effort in \cite{AMP}).

\section{Proof of the Theorem} 

\noindent For the immediately following see the introduction
to \cite{141}, alternatively \cite {Schm}. We recall that if
$x_h/y_h$ is a convergent to
$F$ then 
\begin{equation} \label{eq:conv}
\deg(y_hF-x_h)=-\deg y_{h+1}=-\deg(a_{h+1}y_h+y_{h-1})=-\deg
a_{h+1}-\deg y_h 
\end{equation}
and conversely (Proposition 1 of \cite{141}) that if $x$ and $y$ are
coprime \poly s then
\begin{equation} \label{eq:suff}
\deg(yF-x)\lt-\deg y
\end{equation}
entails $x/y$ is a convergent to $F$.

To prove the Theorem we need to notice that, even though both continuants
$\ov x_h$ and $\ov y_h$ may blow up,  the existence of $\ov
F$ implies that there is some constant $c_h$ in $\K$ so that both
$\ov{c_hx_h}$ and
$\ov{c_hy_h}$ exist and not both vanish identically. 
Then 
$$\deg(\ov{c_hy_h}\ov F-\ov{c_hx_h}\,)\le \deg (y_hF-x_h)\lt -\deg
y_h\le-\deg
\ov{c_hy_h}
$$
so if $r_h=x_h/y_h$ is a convergent to $F$ then, by \eqref{eq:suff}, $\ov
{c_hx_h}/\ov {c_hy_h}=\ov r_h$ is a convergent to
$\ov F$.

Suppose however that the convergents $\ov r_h$, $\ov r_{h+1}$,
$\ldots\,$, $\ov r_{k-1}$ coincide but neither $\ov r_{h-1}=\ov r_h$, nor 
$\ov r_{k-1}=\ov r_{k}$. To show we have found \emph{all} the convergents
to
$\ov F$ it suffices for us to see that $\ov r_h$ and $\ov r_k$ are
\emph{consecutive} convergents to
$\ov F$.

However, of course $\ov r_h-\ov r_k=\ov r_{k-1}-\ov r_k$. 
Note that though neither $x_{k-1}$ and $y_k$, nor $x_k$ and
$y_{k-1}$ might have good reduction, necessarily the difference
$$x_{k-1}y_k-x_ky_{k-1}=(-1)^k$$ and so also the product $y_{k-1}y_k$ must
have good reduction without requiring any multiplier. Thus
$$
\ov{x_{k-1}y_k - x_{k} y_{k-1}}=\pm1
.$$ 
Were it otherwise, we would
have $\ov r_{k-1}=\ov r_k$, contrary to hypothesis. Just so then
$$\ov{x_hy_k-x_k y_h}=\pm1\,,
$$
and necessarily $y_hy_k$, has good reduction, proving that there is no
convergent to $\ov F$ between the convergents
$\ov r_h$ and
$\ov r_k$, as required.

\section{Continued Fraction Manipulation}
\noindent The following remarks are intended to show off several lemmata
allowing one to manipulate \cfe s in somewhat surprising ways. 

In particular we will see by manipulation of \cf s that the blowing
up of a \pq\ in a \cfe\ leads to a collapse to higher degree. The point is
that, even if there is a blowup of \pq s under specialisation or reduction,
one does not have to re-expand the function. In principle one can obtain
the new \cfe\ directly from the original \ex.

\subsection{} One sees readily,  mind you with some pain,
that\begin{multline}\label{eq:first}
F^{-1}=\(aX^{-1}+bX^{-2}+cX^{-3}+dX^{-4}+\cdots\cdots\)^{-1}\\ =
\tfrac 1aX \(1+\tfrac ba X^{-1}+\tfrac ca X^{-2}+\tfrac da
X^{-3}+\cdots\;\)^{-1}\\
=\tfrac 1aX \(1-\tfrac ba X^{-1}-\tfrac ca X^{-2} -\tfrac da X^{-3}-\cdots +
\tfrac
{b^2}{a^2}X^{-2}+\tfrac{2bc}{a^2}X^{-3}+\cdots-\tfrac{b^3}{a^3}X^{-3}-\cdots\)\\
=\(\tfrac 1a X-\tfrac{b}{a^2}\)+\tfrac{b^2-ac}{a^3}X^{-1}
\(1-\tfrac{a^2d-2abc+b^3}{a(b^2-ac)}X^{-1}-\cdots\;\)\\
=\CF{\tfrac 1a
X-\tfrac{b}{a^2}\\\tfrac{a^3}{b^2-ac}X+\tfrac{a^2(a^2d-2abc+b^3)}{(b^2-ac)^2}\\\ldots}.
\end{multline}
However, 
\begin{multline}\label{eq:second}
\(F-aX^{-1}\)^{-1}=\(bX^{-2}+cX^{-3}+dX^{-4}+\cdots\cdots\)^{-1}\\
=
\tfrac 1bX^2 \(1+\tfrac cb X^{-1}+\tfrac db
X^{-2}+\cdots\;\)^{-1}\\ 
=           
\tfrac 1bX^2 \(1-\tfrac cb X^{-1}-\tfrac db
X^{-2}-\cdots+\tfrac{c^2}{b^2}X^{-2}+\cdots\;\)^{-1}\\ 
=\CF{\tfrac 1b X^2-\tfrac{c}{b^2}X+\tfrac{c^2-bd}{b^3}\\\ldots}.
\end{multline}
I now show directly that, indeed, the linear \pq s of the first \ex\ collapse to
a \pq\ of higher degree when $a$ vanishes. To that end I subtract $aX^{-1}$ from
the \cf\ \ex\ of $F$ or, rather --- because this turns out to be more
convenient to do --- I add $aX^{-1}$ to that of $F-aX^{-1}$.

Thus, consider the \ex\
$$F=aX^{-1}+\CF{0\\\tfrac 1b X^2-\tfrac{c}{b^2}X+\tfrac{c^2-bd}{b^3}\\\beta}
=\CF{aX^{-1}\\\tfrac 1b X^2-\tfrac{c}{b^2}X+\tfrac{c^2-bd}{b^3}\\\beta}\,.
$$
Of course, actually to suggest $aX^{-1}$ \emph{is} a \pq\ is ``wash
your mouth out'' stuff and we will have to work hard to make up for
the outrage.  We begin with several fairly well known lemmata, and
their little known corollaries.

\begin{lemma}[Multiplication]\label{le:Schmidt}
$$
B\CF{Ca_0\\Ba_1\\Ca_2\\Ba_3\\Ca_4\\\ldots}=C\CF{Ba_0\\Ca_1\\Ba_2
\\Ca_3\\Ba_4\\\ldots}. 
$$
\end{lemma}
\noindent This fact is both obvious and fairly well known. Its
present  felicitous formulation is given by Wolfgang Schmidt
\cite{Schm}.
\begin{lemma}[Negation]\label{le:negation}
\begin{multline*}
\CF{\alpha\\A\\B\\\beta}=\CF{\alpha\\A\\0\\-1\\1\\-1\\0\\-B\\-\beta}
\\=
\CF{\alpha\\A-1\\1\\-B-1\\-\beta},
\end{multline*}
and
\begin{multline*}
\CF{\alpha\\A\\B\\\beta}=\CF{\alpha\\A\\0\\1\\-1\\1\\0\\-B\\-\beta}
\\=
\CF{\alpha\\A+1\\-1\\-B+1\\-\beta}.
\end{multline*}
\end{lemma}
\begin{proof}The first is just the \ex:
\begin{multline*}
 \CF{-\gamma}= \CF{0\\-1/\gamma}
= \CF{0\\-1\\\gamma/(\gamma-1)}
\\=\CF{0\\-1\\1\\ \gamma-1}
=\CF{0\\-1\\1\\-1\\1/\gamma}
=\CF{0\\-1\\1\\-1\\0\\ \gamma};\end{multline*}
now multiply by $-1$ to get the second claim.
\end{proof}
\begin{corollary}
$$\CF{A+x\\\beta}=\CF{A\\1/x\\-x^2\beta -x} \quad\text{and}\quad
\CF{A\\x\\\beta}=\CF{A+1/x\\-x^2\beta -x}.$$
\end{corollary}
\begin{proof} Just note that
\begin{equation*}
 \CF{x+A\\\beta}= x\CF{1+A/x\\x\beta}
= x\CF{A/x\\1\\-x\beta-1}
=\CF{A\\1/x\\-x^2\beta-x},\end{equation*}
and similarly for the second claim.\end{proof}

\subsection{} We now return to the calculation of the \cfe\ for
$F$. Sequentially we get
\begin{multline*}
F=\CF{aX^{-1}\\\tfrac 1b X^2-\tfrac{c}{b^2}X+\tfrac{c^2-bd}{b^3}\\\beta}
\\=\CF{aX^{-1}\\\tfrac 1b X^2-\tfrac{c}{b^2}X\\\tfrac{b^3}{c^2-bd}\\
-\tfrac{(c^2-bd)^2}{b^6}\beta-\tfrac{c^2-bd}{b^3}}\\
=aX^{-1}\CF{1\\\tfrac a{b} X-\tfrac{ac}{b^2}\\\tfrac{b^3}{a(c^2-bd)}X\\
a\beta'/X}
\\=aX^{-1}\CF{0\\1\\-\tfrac a{b} X+\tfrac{ac}{b^2}-1\\-\tfrac{b^3}{a(c^2-bd)}X\\
-a\beta'/X}
\\=aX^{-1}\CF{0\\1\\-\tfrac a{b}X\\\tfrac{b^2}{ac-b^2}\\
\tfrac{(ac-b^2)^2}{b^4}\tfrac{b^3}{a(c^2-bd)}X-\tfrac{ac-b^2}{b^2}\\
\tfrac{ab^4}{(ac-b^2)^2}\beta'/X}\\
=aX^{-1}\CF{0\\1\\-\tfrac a{b}X\\\tfrac{b^2}{ac-b^2}\\
\tfrac{(ac-b^2)^2}{b^4}\tfrac{b^3}{a(c^2-bd)}X\\
-\tfrac{b^2}{ac-b^2}\\-a\beta'/X+\tfrac{ac-b^2}{b^2}}\\
=\CF{0\\\tfrac1a X\\-\tfrac {a^2}{b}\\\tfrac{b^2}{a(ac-b^2)}X\\
\tfrac{(ac-b^2)^2}{b(c^2-bd)}\\
-\tfrac{b^2}{a(ac-b^2)}X\\-a^2\beta'/X^2+\tfrac{a(ac-b^2)}{b^2}/X}\\
=\CF{0\\\tfrac1a X-\tfrac {b}{a^2}\\-\tfrac{a^3}{ac-b^2}X
+\tfrac {a^2}{b}\\
-\tfrac{b(ac-b^2)^2}{a^4(c^2-bd)}\\
\tfrac{a^3}{ac-b^2}X\\\tfrac{b^2}{a^2}\beta'/X^2-\tfrac{ac-b^2}{a^3}/X}\\
=\CF{0\\\tfrac1a X-\tfrac {b}{a^2}\\-\tfrac{a^3}{ac-b^2}X
+\tfrac {a^2}{b}-\tfrac{a^4(c^2-bd)}{b(ac-b^2)^2}\\
\ldots},
\end{multline*}
which, \emph{mirabile dictu}, is as we had expected and had asserted. 

Note that
the second \ex\ \eqref{eq:second} contains more information than the first
\eqref{eq:first}, whence my reluctance to work from the first \ex. Note
also that the computation just presented is  not the sort of thing one will
essay more than once in a lifetime. Indeed, other than over the finite
field $\F_2$, and perhaps
$\F_3$, it surely cannot be reasonable to attempt to obtain detailed
information in this way in general circumstances. Over those finite fields,
on the other hand, the methods just illustrated may be pursued without pain
or fear; for example, there is work of Niederreiter and Vielhaber
\cite{Nie, Nie1} applying just these notions. Note also the pursuit of
the present ideas over
$\Z$ in a very special case
\cite{134}, and mention of similar such cases in \cite{141}.

\section{Applications}

\subsection{} Denote by $D(X)$ a monic \poly\ of even degree $2g+2$ defined
over $\Z$. Contrary to the numerical case, the \cfe\ of $\sqrt D$ is
\emph{not} periodic in general. Indeed, one proves the existence of a
unit in $\Q(\sqrt D)$, equivalently periodicity of the \ex, by the box
principle. But there are infinitely many \poly s of bounded degree with
coefficients in an infinite field. It is easy to see that the \pq s
$a_h$ of $\sqrt D$ satisfy $\deg a_h\le g$ unless the \ex\ happens to be
periodic in which case the occurrence of a \pq\ of degree $g+1$ signals
the end of a quasi-period. 

For example (see \cite{145}) it happens that
$$\displaylines{
\sqrt{X^4-2X^3+3X^2+2X+1}=
\hfill\cr\hfill
[X^2 - X + 1\,,\,
\overline{\tfrac{1}{2}X - \tfrac{1}{2}\,,\,
2X - 2\,,\,\tfrac{1}{2}X^2 - \tfrac{1}{2}X + \tfrac{1}{2}\,,\,2X - 2\,,\,
\tfrac{1}{2}X -
\tfrac{1}{2}\,,\, 2X^2 - 2X + 2\,}]
}$$
is periodic whereas 
$$
\displaylines{\sqrt{D}=\sqrt{X^4-2X^3+3X^2+2X+2}\hfill
\cr
=
[X^2 - X + 1\,,\,\tfrac{1}{2}X - \tfrac{5}{2^3}\,,\,
\tfrac{2^5}{3\cdot 7}X - \tfrac{2^3\cdot 43}{3^2\cdot 7^2}\,,\,
-\tfrac{3^3\cdot 7^3}{2^8\cdot 31}X - \tfrac{3^2\cdot 7^2\cdot 6719}{2^{11}
\cdot 31^2}\,,\,
\cr-\tfrac{2^{14}\cdot 31^3}{3^4\cdot 7^4 \cdot 13229}X + 
\tfrac{2^{11}\cdot 5^2\cdot 31^2\cdot 329591}{3^4\cdot 7^4 \cdot 13229^2}\,,\,
-\tfrac{3^3\cdot 7^3 \cdot 13229^3}{2^{17}\cdot 5\cdot 31^4 \cdot 1877}X + 
\tfrac{3^2\cdot 7^2 \cdot 13229^2\cdot 21577726507}{2^{19}\cdot 5^2\cdot 31^4
\cdot 1877^2}\,,\,
\cr
-\tfrac{2^{21}\cdot 5^3\cdot 31^4 \cdot 1877^3}{3\cdot 7\cdot 11 
\cdot 13229^4\cdot 12524251}X - \tfrac{2^{19}\cdot 5^2\cdot 31^4\cdot 47\cdot
1877^2\cdot 2693\cdot 1180897}{3^2\cdot 7^2\cdot 11^2\cdot 13229^4\cdot
12524251^2}\,,\,
\hfill\cr
+\tfrac{11^3\cdot 13229^4\cdot 12524251^3}{2^{25}\cdot 5^4\cdot 31^5
\cdot 1877^4\cdot 130960463}X
- \tfrac{11^2\cdot 13229^4\cdot 2109269\cdot 12524251\cdot 208276252871}{2^{29}
\cdot 5^3\cdot 31^6\cdot 1877^4\cdot 130960463^2}\,,\,
\cr
\tfrac{2^{33}\cdot 5^4\cdot 31^7\cdot 1877^4\cdot 130960463^3}{3^2\cdot 7
\cdot 11^4\cdot 67\cdot 331\cdot 13229^4\cdot 12524251^4\cdot 32646599}X 
\hfill\cr\hfill
- \tfrac{2^{29}\cdot 5^4\cdot 31^6\cdot 1877^4\cdot 130960463^2
\cdot 672668401\cdot 6280895711017969}{3^4\cdot 7^2\cdot 11^4\cdot 67^2\cdot
331^2\cdot 13229^4\cdot 12524251^4\cdot 32646599^2}\,,\,
\cr
\hfill \ldots\ldots\;]\,}
$$
seemingly is not.

On the other hand, if we were to view
$D$ as defined over some finite field, say $\F_p$ (with $p\ne2$), then
necessarily its \ex\ must be periodic because now there are only finitely
many \poly s of bounded degree. Indeed, inspection of the
\ex\ above shows that the \pq\ $a_2$ blows up at $p=3$ and at $p=7$, $a_5$
blows up at
$p=5$, $a_6$ at $p=11$, \dots\
 signalling period length $r=2$ at $p=3$ and
$p=7$ (more to the point: regulator, the sum of the degrees of the \pq s
making up a quasi-period, $m=3$ at those primes), regulator $m=6$ at
$p=5$, regulator $m=7$ at $p=11$, \dots\,. This more than suffices,
by a remark of Jing Yu \cite{Yu}, to \emph{prove} that the \ex\ of $\sqrt
{D(X)}$ over
$\Q$ is indeed not periodic.

Specifically, Yu points out that, by the reduction theory of abelian
varieties, if a divisor class (here, that of the divisor at infinity)
is of order $m$ on the Jacobian $\Jac \mathcal C$ of some curve $\mathcal C$
then unless $p\Div m$, it is also  of order $m_p=m$ on $\Jac \mathcal C_p$,
the Jacobian of the curve $\mathcal C$ reduced modulo a prime $p$ of good
reduction for $\mathcal C$; if
$p$ divides
$m$ then $m=m_pp^i$, some positive integer~$i$. For our example
$m_7=3$ and
$m_5=6$ suffices to prove non-periodicity. We should here note that the
primes dividing the discriminant of the \poly\ $X^4-2X^3+3X^2+2X+2$ are
$2$, $3$, and $31$ so both
$7$ and $5$ are primes of good reduction.

\subsection{} Suppose one hopes to find all quartics whose square root
\emph{does} have a periodic \ex. Without too much loss, denote the general
quartic
$D$ by
$$D(X)=(X^2+u)^2+4v(X+w);
$$
here one may also usually suppose the normalisation $u+w^2=v$. It turns out
not to be impossibly painful \cite{163} to compute the \pq s
$a_0(X)=X^2+u$,
$\ldots\,$, 
$a_h(X)=2(X-c_h)/b_h$, by
\begin{equation}\label{eq:b} b_{2h}=\frac{s_3s_5\cdots s_{2h-1}}{s_2s_4s_6\cdots
s_{2h}}
\quad\text{and}\quad
b_{2h+1}=4v\frac{s_2s_4s_6\cdots s_{2h}}{s_3s_5\cdots s_{2h+1}}
\end{equation}
and
\begin{equation}\label{eq:c}
c_{h+1}=(-1)^h\(w-s_2+s_3-s_4+\cdots+(-1)^h s_{h+1}\),
\end{equation}
where
$$s_{h+1}=v/s_h(s_h-1)s_{h-1} \quad\text{for $h=3$, $4$, $\ldots\,$,}
$$
and $s_0=0$, $s_1=\infty$, $s_2=1$, $s_3=v/(1-2w)$. 

One finds that the specialisation $b_{m-1}=0$, equivalently
$s_{m-1}=\infty$, signals that
$a_{m-1}(X)$ blows up (to degree~$2$),  implying regulator $m$. In the
case of periodicity, the `symmetry' $s_1=s_{m-1}$ in fact implies that
$s_2=1=s_{m-2}$,
$s_3=s_{m-3}$, $\ldots\,$. 

The case $m=2$ is $v(X+w)=\kappa$, some constant; $m=3$ is $u=-w^2$.
From there on $u+w^2=v$. Then $m=4$ is $s_3=\infty$ so $1=2w$; $m=5$ is
$s_3=1$ or $v=1-2w$; $m=6$ is $s_4=v/s_3(s_3-1)=1$; $m=7$ is $s_4=s_3$ or
$v=s_3^2(s_3-1)$; $m=8$ is $s_5=s_3$, where in fact
$s_5=(s_3-1)/(s_4-1)$, and so on. 

Indeed, one obtain a raw form of the
modular equation defining $X_1(m)$. For example the condition $s_6=s_5$
for $m=11$ quickly simplifies to 
$$
s_3(s_3-1)=s_5^2(s_5-1).
$$ 
See \cite{163} for details or, for different methods to construct the
families of curves, \cite{Ni},
\cite{AZ}, or
\cite{158}; the original tabulation occurs in \cite{Ku}, see also
\cite{HLP}.

\subsection{} The case $g=2$, thus the issue of the periodicity of
the \cfe\ of the square root of a sextic defined over $\Q$, is 
more interesting if only because, unlike the elliptic case, we do not yet
know all possibilities. Preliminary study suggests it is useful to
distinguish the cases $D(X)$ given by
$$
(X^3+fX+g)^2+u(X^2+vX+w) \quad\text{or}\quad (X^3+fX+g)^2+v(X+w).
$$ 
Generically, after the first several quotients, all \pq s will be of
degree one so, just as in the $g=1$ case, periodicity requires
specialisations of --- that is, relations on --- the coefficients $f$,
$g$, $u$, $v$, $w$. However, now a `blowup' of a \pq\ may mean
no more than that a multiple of the divisor at infinity on the Jacobian
of the hyperelliptic curve $Y^2=D(X)$ is (technically, corresponds to) a
point on the curve and, unlike the elliptic case, does not guarantee
periodicity. It will be interesting to see whether the Reduction Theorem is
helpful in practice in recognising blowup to degree~$g+1$.

\bibliographystyle{amsalpha}

\begin{thebibliography}{ABC}

 




\bibitem{AMP} J.-P. Allouche. M. Mend\`es France, and A. J. van der
Poorten, `An infinite product with bounded partial
quotients\rq, {\it Acta Arith.\/} {\bf 59\/} (1991), 171--182.

\bibitem{AZ}{Roberto M. Avanzi and Umberto M. Zannier, `Genus one curves
defined by separated variable polynomials and a polynomial Pell equation', {\it
Acta Arith.\/} {\bf 99\/}.3 (2001), 227--256.}






\bibitem{Ca}{David G. Cantor, `On the continued fractions of quadratic surds',
{\it Acta Arith.\/} {\bf 68\/} (1994), 295--305.}

\bibitem{HLP} Everett W. Howe, Franck Lepr\'evost, and Bjorn Poonen, `Large
torsion subgroups of split Jacobians of curves of genus two or
three', {\it Forum Math.\/} {\bf 12\/}.3 (2000), 315--364
(MR2001e:11071). 

\bibitem{Ku} Daniel Sion Kubert, `Universal bounds on the torsion of elliptic
curves',  {\it Proc. London Math. Soc.\/} {\bf 33\/}.3 (1976), 193--237. 


\bibitem{Ni} {Abderrahmane Nitaj,  `D\'etermination de courbes
elliptiques pour la conjecture de Szpiro', {\it Acta Arith.\/} {\bf 85\/}.4
(1998), 351--376; `Isogenes des courbes elliptiques definies
sur les rationnels', to appear in {\it J. Combinatorial Math.\/}; see
\url{http://www.math.unicaen.fr/~nitaj/}.}


\bibitem{134} 
Michel Mend\`es France, Alfred J. van der Poorten, and Jeffrey
Shallit, `On lacunary formal power series 
and their continued fraction expansion', in {\it Number Theory in Progress\/} 
(Proc. Conf. in honour of Andrzej Schinzel on the occasion
of his 60th birthday),  K. Gy\H ory, H. Iwaniec and J. Urbanowicz eds, Walter de Gruyter,
Berlin 1999, 321--326.

\bibitem{Nie} Harald Niederreiter and Michael Vielhaber, `Linear complexity
profiles: Hausdorff dimensions for almost perfect profiles and measures for
general profiles',  {\it J. Complexity\/} {\bf 13\/}.3 (1997), 353--383.

\bibitem{Nie1} \bysame, `An algorithm for shifted continued fraction
expansions in parallel linear time. Cryptography', {\it  Theoret. Comput.
Sci.\/} {\bf 226\/}.1--2  (1999), 93--104.


\bibitem{141} {Alf van der Poorten, \lq Formal power series and their
continued fraction expansion', in Joe Buhler ed., {\it Algorithmic Number
Theory\/} (Proc. Third International Symposium, ANTS-III, Portland, Oregon,
June 1998), Springer Lecture Notes in Computer Science {\bf 1423\/} (1998),
358--371.}


\bibitem{145} Alfred J. van der Poorten, 
`Non-periodic
continued fractions in hyperelliptic function fields' (Dedicated to
George Szekeres on his 90th birthday),  {\it Bull. Austral. Math. Soc.\/} {\bf
64\/} (2001), 331--343.

\bibitem{163} \bysame,
`Periodic continued fractions and elliptic curves',  in {\it High
 Primes and Misdemeanours\/}: lectures in honour of the 60th birthday of Hugh
Cowie Williams, Alf van der Poorten and Andreas Stein eds., Fields Institute
Communications {\bf 42\/}, American Mathematical Society, 2004, 353--365. See
also  `Elliptic curves and \cf s', 
\url{www.arxiv.org/math.NT/0403225} for more on this topic.



\bibitem{158} Alfred J. van der Poorten and Xuan Chuong Tran, `Periodic
continued fractions in elliptic function fields', in  Claus Fieker
and David R. Kohel eds, {\it Algorithmic Number Theory\/}  (Proc. Fifth
International Symposium, ANTS-V, Sydney, NSW, Australia July 2002),
Springer Lecture Notes in Computer Science 2369 (2002), 390--404.





\bibitem{Schm} Wolfgang M, Schmidt, `On continued fractions and
diophantine approximation in power series fields', {\it Acta
Arith.\/} {\bf 95\/} (2000), 139--166.


\bibitem{Yu} Jing Yu, `Arithmetic of hyperelliptic curves', manuscript marked
{\it Aspects of Mathematics\/}, Hong Kong University, 1999; see pp.4--6.


\end{thebibliography}


\label{page:lastpage}
\end{document}